\documentclass[12pt, oneside, a4paper]{article}

\usepackage{amsmath}
\usepackage{amsfonts}
\usepackage{amssymb}
\usepackage{amsthm,mathrsfs}
\usepackage{enumerate}
\usepackage{graphicx}
\newtheorem{theorem}{Theorem}[section]
\newtheorem{lemma}[theorem]{Lemma}

\theoremstyle{definition}

\title{\textbf{The Higman-Sims graph and a geometric construction of the Steiner system $S(3,6,22)$}}
\author{Mahdi Ebrahimi\footnote{ m.ebrahimi.math@ipm.ir\\Acknowledgements:
	This research was supported in part
	by a grant from School of Mathematics, Institute for Research in Fundamental Sciences (IPM).}\\ \em ORCID ID: 0000-0001-9789-7376
 \\
 {\small\em  School of Mathematics, Institute for Research in Fundamental Sciences (IPM)},\\{\small\em P.O. Box: 19395--5746, Tehran, Iran}}

\date{}

\begin{document}

\maketitle

\begin{abstract}
The Higman-Sims graph is the unique strongly regular graph with parameters $(100,22,0,6)$. In this paper, we show that the only distance-regular antipodal double cover of the Higman-Sims graph is its bipartite double cover. Our main tool is a geometric construction of the Steiner system $S(3,6,22)$.
\end{abstract}

\section{Introduction}
$\noindent$Association schemes \cite{2} were first defined by Bose and Mesner \cite{8}.
Delsarte used association schemes to unify the study of coding theory and design theory in his thesis \cite{D}.
He identified two types of association schemes, so-called $P$-polynomial and $Q$-polynomial schemes. \textit{distance-regular graphs} are precisely $P$-polynomial association schemes and are well studied.
In particular, imprimitive $P$-polynomial association schemes are distance-regular graphs in which they are bipartite or antipodal \cite[Theorem 4.2.1]{Br}.
A well-known conjecture \cite[p.312]{2} of Bannai and Ito states that for sufficiently large $d$, a primitive $d$-scheme is $P$-polynomial if and only if it is $Q$-polynomial. The study of primitive (imprimitive) distance-regular graphs is a fundamental step towards solving this conjecture. In this paper, we are interested in the classification of distance-regular antipodal covers of strongly regular graphs which are imprimitive $P$-polynomial association schemes.
A \textit{strongly regular graph }can be viewed as an association scheme with two classes. In other words, a strongly regular graph is a distance-regular graph with diameter two.

In this paper, all graphs are assumed to be finite and simple. An important feature of quantum networks is perfect state transfer. Godsil \cite{G} obtained that if a distance-regular graph $\Gamma$ has perfect state transfer, then $\Gamma$ is an antipodal double cover. The Wells graph and Hadamard graphs are examples of distance-regular antipodal double covers. Filo \cite{F} has obtained tow spectral characterizations of this family. Also in \cite{A}, Araya et al. have shown that if $\Gamma$ is a distance-regular graph of valency $k>2$ and diameter $d$ with the parameter $b_t=1$ for $t\leq d/2$ then $\Gamma$ is an antipodal double cover.

 It is well-known that the bipartite double cover of the Gewirtz graph is the unique distance-regular antipodal cover of it \cite[Proposition 11.4.5]{Br}. Nakagawa \cite{M} proved that a distance-regular antipodal double cover with diameter $4$  of a strongly regular graph with parameters $(\nu,k,0,2)$ can be constructed from a symmetric $6$-association scheme. In \cite{Mi}, Miklavic has obtained an upper bound for the valency of distance-regular antipodal covers with diameter $4$ of non-bipartite strongly regular graphs. Also, distance-regular antipodal covers of Hermitian-form graphs have been studied by Makhnev in \cite{M2}. The aim of this paper is to obtain a complete characterization of the distance-regular antipodal double covers of the Higman-Sims graph, which is the unique strongly regular graph with parameters $(100,22,0,6)$.

 \begin{theorem}\label{main}
The only distance-regular antipodal double cover of the Higman-Sims graph is its bipartite double cover.
 \end{theorem}
\section{Distance-regular graphs}
$\noindent$In this section, we present the basic definitions and some known results concerning distance-regular graphs. 
For more background information, we refer the reader to \cite{Br}.
Throughout this section, let $\Gamma$ denote a connected graph with vertex set $V(\Gamma)$. For adjacent vertices $x,y\in V(\Gamma)$, we denote by $xy$, the edge between them.
The \textit{distance} between vertices $x$ and $y$, denoted by $d(x,y)$, is the length of a shortest path between $x$ and $y$ in $\Gamma$.
The \textit{diameter} $\mathrm{diam}(\Gamma)$ of $\Gamma$ is the maximum value of $d(x,y)$ for all pairs of vertices $x$ and $y$ of $\Gamma$. 
Assume that $\mathrm{diam}(\Gamma)=D$.
For $x \in V (\Gamma)$ and an integer $0 \leq i \leq D$, set
$\Gamma_i(x) := \{y \in V (\Gamma) :\, d(x, y) = i\}$. 
The graph $\Gamma$ is called \textit{distance-regular},
if for all integers $0 \leq h, i, j\leq D$ and for all vertices $x, y \in V (\Gamma)$ with $d(x, y) = h$,
the number $p^{h}_{i,j} = |\Gamma_i(x)\cap \Gamma_j(y)|$ is independent of $x$ and $y$.
The numbers $p^{h}_{i,j}$ are called the \textit{intersection numbers} of $\Gamma$. 
Define 
$$c_i:= p^{i}_{1,i-1}, \,\,\,a_i := p^i_{1,i}, \,\,\,b_i := p^i_{1,i+1}, \,\,\,(0\leq i\leq D).$$
$\noindent$Observe that the graph $\Gamma$ is regular with valency $k := b_0$.
Moreover, note that $a_0 = b_D =c_0 = 0, \;c_1 = 1$, and $k=a_i + b_i + c_i$ for every $0 \leq i \leq D$.
The sequence $l(\Gamma):=\{b_0, b_1, \dots , b_{D-1}; c_1, c_2, \dots ,c_D\}$ is referred  as the \textit{intersection array} of $\Gamma$. 
  
   Let $\Gamma$ be distance-regular with diameter $D > 1$. The graph $\Gamma$  is said to be \textit{antipodal}, if the (antipodality) relation "to be at distance $0$ or $D$" is an equivalence relation on $V(\Gamma)$, and the classes of this relation are called \textit{antipodal classes of $\Gamma$}. When $\Gamma$ is antipodal, the \textit{folded graph }$\bar{\Gamma}$ of $\Gamma$ is defined as follows: The graph $\bar{\Gamma}$ has the antipodal classes of $\Gamma$ as vertices, and two equivalence classes $C_1$ and $C_2$ are adjacent whenever the edge set of the induced subgraph of $\Gamma$ on $C_1\cup C_2$ is non-empty. Note that the antipodal classes of $\Gamma$ have the same size $r$ and $\Gamma$ is called an \textit{antipodal $r$-cover} of $\bar{\Gamma}$.
     In addition, when $r=2$, the graph $\Gamma$ is said to be an \textit{antipodal double cover} of $\bar{\Gamma}$. Now, assume that $l(\Gamma)=\{b_0, b_1, \dots , b_{D-1}; c_1, c_2,\dots , c_D\}$. It is well-known that the distance-regular graph $\Gamma$ is an antipodal double cover if and only if for every integer $0\leq i\leq D$, we have $b_i=c_{D-i}$ (see \cite[Proposition 4.2.2]{Br}). If $\tilde{\Gamma}$ is a distance-regular antipodal $r$-cover of $\Gamma$, then the graph $\tilde{\Gamma}$ has diameter $2D$ or $2D+1$ (see \cite[Sect. 4.2]{Br}). Furthermore, if $\mathrm{diam}(\tilde{\Gamma})=2D$, then the graph $\tilde{\Gamma}$ has the intersection array 
  $$\{b_0,b_1,\dots,b_{D-1}, \frac{r-1}{r}c_D,c_{D-1},\dots,c_1;c_1,\dots,c_{D-1}, \frac{1}{r}c_D,b_{D-1},\dots,b_0\}.$$

Let $\Gamma$ be a $k$-regular graph with $\nu$ vertices. The graph $\Gamma$ is said to be \textit{ strongly regular} with parameters $(\nu, k, \lambda, \mu)$, whenever the following conditions hold:
\begin{itemize}
    \item [(i)] $\mathrm{diam}(\Gamma)=2$.
    \item [(ii)] Every pair of adjacent vertices has precisely $\lambda$ common neighbors.
    \item[(iii)] Every pair of nonadjacent vertices has precisely $\mu$ common neighbors.
    \end{itemize}
The\textit{ Higman-Sims graph} is the unique strongly regular graph with parameters $(100,22,0,6)$. Its spectrum is $ \{22^1, 2^{77},(-8)^{22}\}$.  
    This graph was found by Higman and Sims \cite{425}, and uniqueness was proved by Gewirtz \cite{341}.

    The Higman-Sims graph can be constructed as follows: Take a symbol $\infty$, the $22$ points, and the $77$ blocks of the Steiner system $S(3,6,22)$ (see Sect. 4 below) as the $1+22+77=100$ vertices. The vertex $\infty$ is adjacent to the points, a point is adjacent to the blocks containing it, and two blocks are adjacent when they are disjoint.     
    For details, we refer the reader to \cite[Sect. 10.31]{strongly}.
    
    \begin{lemma}\cite[Sect. 4.2]{Br}\label{double}.
 Assume that $\Gamma$ is a strongly regular graph with parameters $(\nu, k, \lambda,\mu)$. If $\tilde{\Gamma}$ is a distance-regular antipodal double cover of $\Gamma$ with $\mathrm{diam}(\tilde{\Gamma})=4$, then 
 \begin{itemize}
     \item [(i)] The eigenvalues of $\Gamma$ are the eigenvalues (old eigenvalues) of $\tilde{\Gamma}$ with the same multiplicities. 
     \item [(ii)] $\tilde{\Gamma}$ has two new eigenvalues $\theta_1$ and $\theta_2$ in which they are the two roots of $\theta^2-\lambda\theta-k=0$, and occur with multiplicity $m(\theta)=\frac{\nu}{2+\lambda\theta/k}$, where $\theta\in \{\theta_1,\theta_2\}$. Consequently, either $\lambda=0$, or $\lambda^2+4k$ is a square (and these eigenvalues are integers).
 \end{itemize}
 \end{lemma}
   
 Let $V(\Gamma)=\{v_1,v_2,\dots,v_n\}$. The \textit{bipartite double cover} of the graph $\Gamma$ is a graph in which it has two vertices $u_i$ and $w_i$ for each vertex $v_i$ of $\Gamma$. Two vertices $u_i$ and $w_j$ are joined by an edge in the bipartite double cover if and only if $v_i$ and $v_j$ are adjacent in $\Gamma$. 
 \begin{lemma}\label{bdc}\cite{M}.
     Let $\tilde{\Gamma}$ be a distance-regular antipodal double cover of the Higman-Sims graph. If $\mathrm{diam}(\tilde{\Gamma})=5$, then $\tilde{\Gamma}$ is the bipartite double cover of the Higman-Sims graph.
 \end{lemma}
\section{ A geometric construction of the Steiner system $S(3,6,22)$}
In this section, we present a geometric construction of the Steiner system $S(3,6,22)$. We require some definitions and notations from \cite{strongly}. A \textit{$t-(\nu, k, \lambda)$ design} is a set of $\nu$ points together with a collection of subsets of size $k$, called\textit{ blocks}, such that each set of $t$ points is precisely in $\lambda$ blocks. A \textit{Steiner system} $S(t,k,\nu)$ is such a design with $\lambda=1$. There exist unique Steiner systems $S(5,8,24)$, $S(4,7,23)$ and $S(3,6,22)$ which are called \textit{Witt designs}. The \textit{derived design} of a $t-(\nu,k,\lambda)$ design $\mathcal{D}$ is a design obtained by deleting one point and all blocks not containing that point from $\mathcal{D}$ which is a $(t-1)-(\nu-1,k-1,\lambda)$ design. 

Let $F$ be a field and $V$ be an $(n+1)$-dimensional vector space over $F$. The set of all subspaces of $V$ is called the $n$-dimensional \textit{projective space} over $F$. In addition, the $n$-dimensional projective space over the finite field $\mathrm{GF}(q)$ is denoted by $\mathrm{PG}(n,q)$. The points, lines and planes of $\mathrm{PG}(n,q)$ are the $1$-dimensional, $2$-dimensional and $3$-dimensional subspaces of $V$, respectively.

Since the Steiner system $S(3,6,22)$ is the derived design of $S(4,7,23)$, for our construction, we closely follow \cite{Ba}. Let $S:=\{1,2,3,\dots,8\}$. Fix a Steiner system $S(3,4,8)$ on $S$, say $D_0$. For every permutation $\pi$ on $S$, we will denote by $D^{\pi}_0$, the image of $D_0$ under $\pi$. Assume that $\mathcal{D}$ is the collection of the $15$ distinct Steiner systems $D_0$, $D_0^{(a,b,c)}$ and $D_0^{(a,c,b)}$, where $B=\{1,a,b,c\}$ ranges over blocks of $D_0$ containing $1$. 
\begin{lemma}\label{block}\cite[Theorem 2.4]{Ba}.
	If $B$ is a block of one of the elements of $\mathcal{D}$  then $B$ is a block of precisely three elements of $\mathcal{D}$. Moreover, the intersection of any two of the elements of $\mathcal{D}$ has size $2$ and contains a block and its complement. There is a unique third element of $\mathcal{D}$ containing this intersection.
\end{lemma}
Define a binary operation $+$ on $\mathcal{D}\cup \{0\}$ as follows:
\begin{itemize}
	\item [(i)] For every $D\in \mathcal{D}$, define $D+D:=0$ and $D+0:=D$.
	\item [(ii)] For each distinct $D$ and $D^{'}$ in $\mathcal{D}$, the Steiner system $D+D^{'}$ is the unique third element of $\mathcal{D}$ containing $D\cap D^{'}$ described in Lemma \ref{block}.
\end{itemize}
$\noindent$ Note that $(\mathcal{D}\cup \{0\},+)$ is a $4$-dimensional vector space over $GF(2)$ \cite[Proposition 2.7]{Ba}.

For the construction of the Steiner system $S(4,7,23)$, take as the point set, the $15$ points of $\mathrm{PG}(3,2)$ together with the elements of $S$. Refer to the $15$ points of $\mathrm{PG}(3,2)$ as \textit{geometric points} and the $8$ points of $S$ as \textit{distinguished points}. The blocks are of three types. These blocks are $15$ planes of $\mathrm{PG}(3,2)$ together with $70$ line sets and $168$ extended ovoid sets, which are defined as follows. We denote by $\mathcal{D}^*$, the $3$-dimensional projective space over $\mathrm{GF}(2)$ defined on $\mathcal{D}\cup \{0\}$. Let $\phi:\mathrm{PG}(3,2)\to \mathcal{D^*}$ be an isomorphism of projective spaces.
Then for each element $D$ of $\mathcal{D}$ and each block $B$ of $D$ assign to $B$ the geometric point $\phi^{-1}(D)$. Since $B$ precisely occurs in three elements of $\mathcal{D}$, the block $B$ is assigned three geometric points which are on a line $L_B$ of $\mathrm{PG}(3,2)$. Define the \textit{line sets} to be all the sets of the form $B\cup L_B$, where $B$ ranges over $70$ blocks of $15$ elements of $\mathcal{D}$. To construct the extended ovoid sets, assume that $a$ and $b$ are distinguished points and $\alpha$ is a geometric point. The triple $(a,b,\alpha)$ occurs exactly three times in the line sets. Let $\pi_\alpha$ be the set of all geometric points of these three line sets. Note that $\pi_\alpha$ is a plane of $\mathrm{PG}(3,2)$ \cite[Proposition 3.1]{Ba}. For distinguished points $a$ and $b$, choose geometric points $\alpha$ and $\beta$ which do not occur together with $a$ and $b$ in the line sets. Define the \textit{extended ovoid sets} to be all the sets of the form $B(a,b,\alpha,\beta):=\{a,b,\alpha,\beta\}\cup (\bar{\pi}_\alpha\cap\bar{\pi}_\beta\backslash\{\alpha+\beta\}) $, where $\bar{\pi}_\alpha$ and $\bar{\pi}_\beta$ are the complements of $\pi_\alpha$ and $\pi_\beta$, respectively. By \cite[Theorem 3.11]{Ba}, the $15$ planes of $\mathrm{PG}(3,2)$ together with $70$ line sets and $168$ extended ovoid sets are precisely the blocks of the Steiner system $S(4,7,23)$. If we delete the distinguished point $1$ from $S(4,7,23)$, then we obtain  the Steiner system $S(3,6,22)$. In particular, if we remove the distinguished point $1$ from line sets containing this distinguished point, we obtain $35$ blocks of $S(3,6,22)$. We call them \textit{linear blocks} of $S(3,6,22)$.  
\section{$2$-lifts and strongly regular signed graphs}
  In this section, we review some definitions and notations from \cite{B} and \cite{stanic}. Let $G$ be a finite simple graph with the vertex set $V$ and the edge set $E$. A \textit{signed graph} $(G,\sigma)$ of $G$ is obtained from the graph $G$ by accompanying each edge $e$ by the sign $\sigma(e)\in\{1,-1\}$. We say that $G$ is the\textit{ underlying graph }of $(G,\sigma)$. The set of edges of $(G,\sigma)$ is composed of subsets of positive and negative edges. Now, let $n:=|V|$. The $n\times n$ \textit{adjacency matrix} $A(G,\sigma)$ of $(G,\sigma)$ is obtained from the standard $(0,1)$-adjacency matrix of $G$ by reversing the sign of all $1$'s that correspond to negative edges. The \textit{$2$-lift} $\mathrm{Lift}(G,s)$ of the graph $G$ associated with a signature $s:E\to \{\pm1\}$ is a graph, defined as follows: The vertices of $\mathrm{Lift}(G,s)$ associated with every vertex $x\in V$ are two vertices, $x_0$ and $x_1$. If $xy\in E$ is a positive edge of $(G,s)$, then the corresponding edges in $\mathrm{Lift}(G,s)$ are $x_0y_0$ and $x_1y_1$. In addition, if $xy$ is a negative edge of $(G,s)$, then the corresponding edges in $\mathrm{Lift}(G,s)$ are $x_0y_1$ and $x_1y_0$. A detailed introduction to $2$-lifts of graphs can be found in \cite{B}.
 \begin{lemma}\cite[Lemma 3.1]{B}.\label{lift}
     Let $A$ be the adjacency matrix of a graph $G$, and $(G,\sigma)$ be a signed graph with the underlying graph $G$. Then every eigenvalue of $A$ and every eigenvalue of $A(G,\sigma)$ are eigenvalues of $\mathrm{Lift}(G, \sigma)$. Furthermore, the multiplicity of each eigenvalue of $\mathrm{Lift}(G, \sigma)$ is the sum of its multiplicities in $A$ and $A(G,\sigma)$.
 \end{lemma}

 A \textit{switching} of a signed graph $(G,\sigma)$ at the vertex $v$ is the change of the signs of all incident edges with $v$. A \textit{switching} for the signed graph $(G,\sigma)$ is a finite sequence of switching at some vertices of $(G,\sigma)$. Two signed graphs $(G_1,\sigma_1)$ and $(G_2,\sigma_2)$ are said to be \textit{switching equivalent}, if $(G_2,\sigma_2)$ is a switching of $(G_1,\sigma_1)$. It can be easily seen that if two signed graphs $(G_1,\sigma_1)$ and $(G_2,\sigma_2)$ with the same underlying graph $G$ are switching equivalent, then the graphs $\mathrm{Lift}(G_1,\sigma_1)$ and $\mathrm{Lift}(G_2,\sigma_2)$ are isomorphic.
 
 We say that a signed graph is \textit{complete}, \textit{totally disconnected} or \textit{bipartite} if the same holds for its underlying graph. A signed graph is called \textit{homogeneous} if all edges have the same sign. A \textit{walk} in a signed graph is a sequence of alternate vertices and edges in which consecutive vertices are the endpoints of the corresponding edge. If the number of negative edges (counted with their multiplicity if there are repeated edges) of a walk in a signed graph is even, then it is called a \textit{positive walk}. Otherwise, it is negative. The difference between the number of positive and negative walks that traverse along $2$ edges between the vertices $x$ and $y$ is indicated by $w_2(x,y)$.

A signed graph $(G,\sigma)$ is said to be \textit{strongly regular} \cite{stanic}, if it is regular and satisfies the following conditions:
\begin{itemize}
\item[(i)]  $(G,\sigma)$ is neither homogeneous complete nor totally disconnected.
\item[(ii)] There exists an integer $a$ such that $w_2(x,y)=a$ for all adjacent vertices $x$ and $y$ with $\sigma(xy)=+1$. 
\item[(iii)] There exists an integer $b$ such that $w_2(x,y)=b$ for all adjacent vertices $x$ and $y$ with $\sigma(xy)=-1$.
\item[(iv)] There exists an integer $c$ such that $w_2(x,y)=c$ for all non-adjacent vertices $x$ and $y$ of $G$.  
\end{itemize}
We end this section with the following result on signed graphs with exactly two distinct eigenvalues. 
\begin{lemma}\label{two}\cite[Theorem 4.1]{stanic}.
An inhomogeneous signed graph has $2$ distinct eigenvalues if and only if it is a strongly regular signed graph with $b=-a$, and if it is non-complete, then $c=0$.
\end{lemma}

\section{Distance-regular antipodal double covers of the Higman-Sims graph}
 The aim of this section is to prove Theorem \ref{main}  our main tool is the geometric construction of the Steiner system $S(3,6,22)$ described in Section 3.\\
\textit{Proof of Theorem \ref{main}:} Assume that $\tilde{\Gamma}$ is a distance-regular antipodal double cover of the Higman-Sims graph $\mathcal{HS}$. Then $\mathrm{diam}(\tilde{\Gamma})=4\,or\, 5$. 
It follows from Lemma \ref{bdc} that we must show that $\mathrm{diam}(\tilde{\Gamma})= 5$. On the contrary, suppose $\mathrm{diam}(\tilde{\Gamma})= 4$. Then as the Higman-Sims graph is strongly regular with parameters $(100,22,0,6)$, the graph $\tilde{\Gamma}$ is a distance-regular graph with intersection array $\{22,21,3,1;1,3,21,22\}$. In addition, by Lemma \ref{double}, we deduce that the spectrum of $\tilde{\Gamma}$ is $\{22^1,(\sqrt{22})^{50}, 2^{77},(-\sqrt{22})^{50},(-8)^{22}\}$. Since $\mathcal{HS}$ is the folded graph of the connected graph $\tilde{\Gamma}$ and $\mathrm{diam}(\tilde{\Gamma})= 4$, there exists an inhomogeneous signed graph $(\mathcal{HS}, \sigma)$ of the Higman-Sims graph $\mathcal{HS}$ such that the graph $\tilde{\Gamma}$ is isomorphic to $\mathrm{Lift}(\mathcal{HS},\sigma)$.
 Hence it follows from Lemma \ref{lift} that the signed graph $(\mathcal{HS}, \sigma)$ has precisely two distinct eigenvalues $\pm\sqrt{22}$. Therefore, using Lemma \ref{two}, the signed graph $(\mathcal{HS}, \sigma)$ is strongly regular with parameters $a=b=c=0$. 

Note that the vertex set of $\mathcal{HS}$ consists of the symbol $\infty$, the $15$ geometric points of $\mathrm{PG}(3,2)$, the $7$ distinguished points in $T:=\{2,3,\dots,8\}$, and the $77$ blocks of the Steiner system $S(3,6,22)$. Since for every two distinguished points $ i,j\in T$, the vertices $i$ and $j$ have $6$ common neighbors, there exist exactly $5$ blocks of $S(3,6,22)$ containing $\{i,j\}$. 
It is easily checked that we must choose  these $5$ blocks from linear blocks described in Section 3. We denote them by $b_{1ij}, b_{2ij},\dots,b_{5ij}$. We know that if two signed graphs $(\Gamma_1,\sigma_1)$ and $(\Gamma_2,\sigma_2)$ with the same underlying graph are switching equivalent, then the graphs $\mathrm{Lift}(\Gamma_1,\sigma_1)$ and $\mathrm{Lift}(\Gamma_2,\sigma_2)$ are isomorphic. Thus, without loss of generality, we can assume that
\begin{itemize}
    \item[(i)] For every distinguished point $i$, we have $\sigma(\infty i)=+1$.
    \item[(ii)] For every liner block $b$, there exist distinct distinguished points $i,j\in b$ such that $\sigma(bi)=\sigma(bj)=+1$. 
\end{itemize}
$\noindent$Then as $(\mathcal{HS},\sigma)$ is a strongly regular signed graph with parameters $a=b=c=0$, for each two distinguished points $i,j\in T$ with $i\neq j$, we get
$$\sum_{t=1}^5\sigma(b_{tij}i)\sigma(b_{tij}j)=-1.$$
$\noindent$Hence for every $i,j\in T$ with $i\neq j$, there exist exactly two linear blocks $b,b'\in \{b_{tij}: 1\leq t\leq 5\}$ such that the edges $bi,b'i,bj$ and $b'j$ are positive edges of the signed graph $(\mathcal{HS},\sigma)$.

By the above discussion, if for every linear block $b$, we define 
$$b^*:=(b\cap T)\backslash \{i\in T:\,\sigma(bi)=-1\},$$
then for every two distinguished points $i,j\in T$ with $i\neq j$, the set $\{i,j\}$ appears twice in the sets $b^*$'s, where $b$ ranges over linear blocks.
Let $p$ (resp. $q$) be the number of linear blocks $b$ satisfying $|b^*|=2$ (resp. $|b^*|=3$). Then as the number of $2$-element subsets of $T$ is equal to $21$, we can see that $p+q=35$ and $p+3q=42$. This is a contradiction. Thus, $\mathrm{diam}(\tilde{\Gamma})=5$ and this completes the proof. \qed
\section*{Declaration of competing interests}
I declare that I have no known competing financial interests or personal relationships that could have appeared to influence the work reported in this paper. 

\section*{Data availability}
Data sharing is not applicable to this article as no new data was created or analyzed in this study.


\end{document}